\documentclass{elsarticle}
\usepackage{amsmath}
\usepackage{amsthm}
\usepackage{amssymb}
\usepackage{graphicx}

\newtheorem{theorem}{Theorem}[section]
\newtheorem{lemma}[theorem]{Lemma}
\newtheorem{prop}[theorem]{Proposition}
\newtheorem{corollary}[theorem]{Corollary}
\newtheorem{defn}[theorem]{Definition}
\newtheorem{ex}[theorem]{Example}

\newcommand{\codim}{\mathrm{codim}}

\begin{document}
\makeatletter
\def\ps@pprintTitle{%
 \let\@oddhead\@empty
 \let\@evenhead\@empty
 \def\@oddfoot{\centerline{\thepage}}%
 \let\@evenfoot\@oddfoot}
\makeatother

\title{Gotzmann regularity for globally generated coherent sheaves}
\author[RD]{Roger Dellaca\corref{cor1}}\ead{rdellaca@uci.edu}
\cortext[cor1]{Corresponding author}
\address[RD]{University of California Irvine, Mathematics Department, Irvine, CA 92697, USA}
\begin{abstract}
In this paper, Gotzmann's Regularity Theorem is established for globally generated coherent sheaves on projective space. This is used to extend Gotzmann's explicit construction to the Quot scheme. The Gotzmann representation is applied to bound the second Chern class of a rank 2 globally generated coherent sheaf in terms of the first Chern class.
\end{abstract}


\maketitle

\section{Introduction}

Gotzmann gave an explicit construction for the Hilbert scheme $\mathrm{Hilb}_P(\mathbb{P}^n)$ of subschemes of $\mathbb{P}^n$ with Hilbert polynomial $P$, in \cite{Got78}. He used a binomial expansion of the Hilbert polynomial to give his Regularity and Persistence theorems as part of this construction; this expansion is called the Gotzmann representation of $P$. The number of terms in the Gotzmann representation is the Gotzmann number of the Hilbert polynomial. This construction gives $\mathrm{Hilb}_P(\mathbb{P}^n)$ as a subscheme of a cross product of Grassmannians which is determined by the Gotzmann number. Alternatively, it can be given as the degeneracy locus of a single Grassmannian determined by the Gotzmann number.  

Some of Gotzmann's results have been extended to other settings; for example, Gasharov \cite{Gas97} proved Gotzmann Persistence for modules. A natural question is: Does Gotzmann Regularity hold for modules? And if so, does Gotzmann's Hilbert scheme construction have an analog for the Quot scheme? Such a construction may have applications to problems in enumerative geometry that reduce to intersections on the Quot scheme, by computing the intersections in the Grassmannian using Schubert calculus; Donaldson-Thomas invariants are a potential example \cite{Mar07}.

The author will show that the Hilbert polynomial of an arbitrary module does not have a Gotzmann representation, but a Gotzmann representation exists for the Hilbert polynomial of a globally generated coherent sheaf, and use this to extend Gotzmann's Regularity Theorem to this class of sheaves. This will be used to extend Gotzmann's construction to the Quot scheme $\mathrm{Quot}_P(\mathcal{O}_{\mathbb{P}^n}^r)$ of quotients of $\mathcal{O}_{\mathbb{P}^n}^r$ with Hilbert polynomial $P$. An example will show that Gotzmann regularity is not a sharp bound on Castelnuovo-Mumford regularity, as it is in the case of subschemes. In addition, the Gotzmann representation will be applied to bound the second Chern class of a rank 2 globally generated coherent sheaf on $\mathbb{P}^3$ in terms of its first Chern class.

\section{Preliminaries}

Let $k$ be a field, $S = k[x_0, \ldots, x_n]$, and $\mathbb{P}^n = \mathrm{Proj}(S)$. Let $F = Se_1 + \cdots + Se_r$ be the free $S$-module of rank $r$ with $\deg(e_i) = f_i$ and $f_1 \leq \cdots \leq f_r \leq 0$. A quotient of $F$ corresponds to a coherent sheaf on $\mathbb{P}^n$ that is generated by global sections.

Let $N$ be a graded $S$-module. The \emph{Hilbert function} of $N$ is $H(N,d) = \dim_k(N_d)$, the dimension of the degree-$d$ part of $N$ as a $k$-vector space. For $d \gg 0$, the Hilbert function becomes a polynomial $P_N(d)$, the \emph{Hilbert polynomial} of $N$. For a coherent sheaf $\mathcal{F}$ on $\mathbb{P}^n$, the Hilbert function and Hilbert polynomial of $\mathcal{F}$ are $H(\mathcal{F},d) = H(\Gamma^*(\mathcal{F}),d)$ and $P_{\mathcal{F}}(d) = P_{\Gamma^*(\mathcal{F})}(d)$.

Given $a, d \in \mathbb{N}$, the \emph{$d$th Macaulay representation of $a$} is the unique expression
\[
a = {k(d) \choose d} + {k(d-1) \choose d-1} + \cdots + {k(\delta) \choose \delta},
\]
with $\delta \in \mathbb{Z}$, satisfying $k(d) > \cdots > k(\delta) \geq \delta > 0$. Given this representation, the \emph{$d$th Macaulay transformation of $a$} is
\[
a^{\langle d \rangle} = {k(d) + 1 \choose d + 1} + {k(d-1) + 1 \choose d} + \cdots + {k(\delta) + 1 \choose \delta + 1}.
\]

\begin{ex} The 3rd Macaulay representation of 11 is
\[
{5 \choose 3} + {2 \choose 2},
\]
and
\begin{align*}
11^{\langle 3 \rangle} & = {6 \choose 4} + {3 \choose 3} \\
   & = 16.
\end{align*}

\end{ex}

\begin{defn}
Given a polynomial $P(d) \in \mathbb{Q}[d]$, a \emph{Gotzmann representation} of $P$ is a binomial expansion
\[
P(d) = {d + a_1 \choose d} + {d+a_2 - 1 \choose d-1} + \cdots + {d + a_s - (s-1) \choose d-(s-1)},
\]
with $a_1, \cdots, a_s \in \mathbb{Z}$ and $a_1 \geq \cdots \geq a_s \geq 0$.
\end{defn}

\begin{prop}[\cite{Vas04}, Corollary B.31]
The Hilbert polynomial of a subscheme of $\mathbb{P}^n$ has a unique Gotzmann representation. $\qed$
\end{prop}

\begin{defn}
The number of terms in the Gotzmann representation of a scheme's Hilbert polynomial is called the \emph{Gotzmann number} of a scheme. 
\end{defn}

\begin{ex} Assume $X$ is a scheme with Hilbert polynomial $P_X(d) = 3d+2$.  This has Gotzmann representation
\[
3d + 2 = {d+1 \choose d} + {d \choose d-1} + {d-1 \choose d-2} + {d-3 \choose d-3} + {d-4 \choose d-4},
\]
and the Gotzmann number of $P_X(d)$ is $5$. 
\end{ex}

Before giving Gotzmann's Regularity Theorem, we must define Castelnuovo-Mumford regularity.

\begin{defn}
A coherent sheaf $\mathcal{F}$ on $\mathbb{P}^n$ is called \emph{$m$-regular} if 

$
H^i(\mathbb{P}^n, \mathcal{F}(m-i)) = 0 \text{ for all } i > 0.
$

The \emph{Castelnuovo-Mumford regularity} (written $\mathrm{reg} \mathcal{F}$) of $\mathcal{F}$ is the smallest such $m$.
\end{defn}

It is known that if $\mathrm{reg} \mathcal{F} = s$, then $\mathcal{F}$ is $d$-regular for all $d \geq s$.

We are ready to give Gotzmann's Regularity and Persistence theorems. 

\begin{theorem}[Gotzmann's Regularity Theorem, \cite{Got78}]
If $X$ is a projective k-scheme with Gotzmann number $s$, then $\mathcal{I}_X$ is $s$-regular. $\qed$
\end{theorem}

\begin{theorem}[Gotzmann's Persistence Theorem, \cite{Got78}]
Suppose $I$ is a graded ideal of $S$ generated in degree at most $d$. If $H(S/I,d+1) = H(S/I,d)^{\langle d \rangle}$, then 

$
H(S/I, s+1) = H(S/I,s)^{\langle s \rangle}
$
for all $s \geq d$. $\qed$
\end{theorem}

Gasharov extended Gotzmann's Persistence Theorem to the case of modules in \cite{Gas97}. 

\begin{theorem}[\cite{Gas97}, Theorem 4.2] \label{Gas97thm}
Assume $F$ is a free graded $S$-module with $l$ the maximal degree of its generators,
$N$ a submodule of $F$ and $M = F /N$. For each pair $(p, d)$ such that $p \geq 0$
and $d \geq p + l + 1$, we have:
\begin{enumerate}
\item $H(M, d + 1) \leq H(M, d)^{\langle d - l - p \rangle}$; 
\item If $N$ is generated in degree at most $d$ and
$H(M, d + 1) = H(M, d)^{\langle d - l - p \rangle}$,
then
$H(M, d + 2) = H(M, d + 1)^{\langle d+1-l-p \rangle}$. $\qed$
\end{enumerate}
\end{theorem}

However, the next example shows that not all Hilbert polynomials of graded $S$-modules have Gotzmann representations.

\begin{ex}
Consider $\mathcal{O}_{\mathbb{P}^1}(-1)$. It has Hilbert polynomial $P(d) = d$. By reason of degree, the Gotzmann representation would have to have ${d+1 \choose 1} = d+1$ as its first term; and the Gotzmann representation can only add more positive terms, so no Gotzmann representation exists.
\end{ex}

\section{Gotzmann representations for globally generated sheaves}

As before, assume a graded module $F = Se_1 + \cdots + Se_r$ with $\deg(e_i) = f_i$ and $f_1 \leq \cdots \leq f_r \leq 0$. This section will establish the following proposition.

\begin{prop} \label{prop1} If $M$ is a submodule of $F$, then $P_{F/M}(d)$ has a unique Gotzmann representation.
\end{prop}

To prove this proposition, some lemmas are necessary.

\begin{theorem}[\cite{Mac27}] Let $H:\mathbb{N} \to \mathbb{N}$ and let $k$ be a field. The following are equivalent. 
\begin{enumerate}
  \item $H$ is the Hilbert function of a projective $k$-scheme.
  \item $H(0) = 1$ and $H(d)^{\langle d \rangle} \leq H(d+1)$ for all $d$. $\qed$
\end{enumerate}
\end{theorem}

\begin{lemma} \label{lemma1} The following are equivalent.
\begin{enumerate} 
  \item $P$ is the Hilbert polynomial of a projective $k$-scheme.
  \item $P$ has a Gotzmann representation.
\end{enumerate}
\end{lemma}

\proof $1 \Rightarrow 2$ is given by Gotzmann's Regularity Theorem. Conversely, if
\[
P(d) = {d+a_1 \choose a_1} + {d+a_2 -1 \choose a_2} + \cdots + {d+a_s -(s-1) \choose a_s},
\]
then define $H(d)$ as
\[
H(d) = {d+a_1 \choose a_1} + {d+a_2 -1 \choose a_2} +  \cdots + {a_{d+1} \choose a_{d+1}}
\]
for $d \leq s-2$, and $H(d) = P(d)$ for $d \geq s-1$. Then 
\[
H(0) = {a_1 \choose a_1} = 1.
\]
If $d \leq s-2$, then
\begin{align*}
H(d+1) & = {d+1+a_1 \choose a_1} + \cdots + {a_{d+1} + 1  \choose a_{d+1}} \\
& \ \ + {d+1+a_{d+2} - (d+1) \choose a_{d+2}} \\
   & = {d+ a_1 +1 \choose d+1} + \cdots + {a_{d+1} + 1 \choose 1} + {a_{d+2} \choose 0} \\
   & = H(d)^{\langle d \rangle} + {a_{d+2} \choose 0} \\
   & \geq H(d)^{\langle d \rangle};
\end{align*}
and if $d \geq s-1$, then 
\[
H(d+1) = H(d)^{\langle d \rangle}.
\]
By Macaulay's criterion, $H$ is the Hilbert function of a projective $k$-scheme, and $H$ has Hilbert polynomial $P$. $\qed$

\begin{lemma} \label{lemma2} If $P$ and $Q$ are polynomials with Gotzmann representations, then $P+Q$ has a Gotzmann representation.
\end{lemma}

\proof By Lemma \ref{lemma1}, $P$ and $Q$ are Hilbert polynomials of projective $k$-schemes $X$ and $Y$ respectively. There exists a projective space $\mathbb{P}^N$ such that $X$ and $Y$ can be embedded disjointly. Then $X \cup Y$ is a projective $k$-scheme with Hilbert polynomial $P+Q$, so by Lemma \ref{lemma1}, $P+Q$ has a Gotzmann representation. $\qed$

\begin{lemma} \label{lemma3} The polynomial 
\[
{d+n+1 \choose n}
\]
 has a Gotzmann representation.
\end{lemma}

\proof One easily sees that 
\[
{d+n+1 \choose n} = 1 + {d+1 \choose 1} + {d+2 \choose 2} + \cdots + {d+n \choose n}.
\]
Lemma \ref{lemma2} then gives the result. $\qed$

\proof (of Proposition \ref{prop1}): By Hulett \cite{Hul95} and Pardue \cite{Par94}, we can assume $F/M = S/I_1 e_1 + \cdots + S/I_r e_r$. There exist Gotzmann representations 
\[
P_{S/I_j}(d) = {d+a_{j,1} \choose a_{j,1}} + {d+a_{j,2} -1 \choose a_{j,2}} + \cdots + {d+a_{j,s_j} -(s_j-1) \choose a_{j,s_j}}.
\]
It is sufficient to show that $P_{S/I_j}(d+1)$ has a Gotzmann representation. Note that
\[
P_{S/I_j}(d+1) = {d+a_{j,1} +1 \choose a_{j,1}} + {d+a_{j,2}  \choose a_{j,2}} + \cdots + {d+a_{j,s_j} -(s_j-2) \choose a_{j,s_j}};
\]
The first term has a Gotzmann representation by Lemma \ref{lemma3}; the remaining terms are a Gotzmann representation; the sum of these two polynomials has a Gotzmann representation by Lemma \ref{lemma2}. By Lemma \ref{lemma1}, this is the Hilbert polynomial of some projective $k$-scheme, and by Proposition \ref{prop1} the Gotzmann representation is unique.
$\qed$

\section{Gotzmann regularity}

We can now prove Gotzmann regularity for globally generated sheaves. The argument is very similar to the one given for the case of ideals in \cite{Vas04}.

\begin{prop} \label{prop2} If $F$ is a rank-$r$ free $S$-module with module generators having degree at most zero, $N$ a graded submodule of $F$ and $M = F/N$, with Gotzmann representation
\begin{equation} \label{eq1}
P_M (d) = {d + a_1 \choose a_1}
+
{d + a_2 - 1 \choose a_2}
+ \cdots +
{d + a_s - (s - 1) \choose a_s},
\end{equation}
then the associated sheaf $\tilde{N}$ is $s$-regular.
\end{prop}

\proof As in \cite{Gas97}, we can assume that the largest degree of a module generator of $F$ is $0$; otherwise, twist $F$ and $N$ by the same amount to cause it, and such twist only makes the Castelnuovo-Mumford regularity smaller by the amount of the twist. We can also assume $N$ is saturated, and that $N \neq F$.

Induct on the number $n= \dim S - 1$. The case $n=0$ is trivial.

In the case $n \geq 1$, there exists $h \in S_1$ such that $h$ is $M$-regular. Indeed, if no $h \in S_1$ is $M$-regular, then $(S_1) = \cup_{h \in S_1}h$ is a zero-divisor on $M$, so there exists a non-zero $z\in F \setminus N$ such that $S_1 z \subset N$, contradicting $N$ saturated.

Set 
\[
S' = S/hS(-1), \ F' = F/hF(-1), \ N' = N/hN(-1), \ M' = M/hM(-1). 
\]

Then $M'$ satisfies the hypotheses of Proposition \ref{prop1}, so we can write
\begin{equation} \label{eq2}
P_{M'} (d) = {d + b_1 \choose b_1}
+
{d + b_2 - 1 \choose b_2}
+ \cdots +
{d + b_r - (r - 1) \choose b_r},
\end{equation}
and $\tilde{N'}$ is $r$-regular.

Also, from the exact sequence $0 \to M(-1) \stackrel{h}{\to} M \to M' \to 0$, it follows that 
\begin{equation} \label{eq3}
P_{M'} (d) = P_M (d) - P_M (d-1).
\end{equation}
Equations (\ref{eq1}), (\ref{eq2}) and (\ref{eq3}) give

\begin{equation} \label{eq4}
P_{M'} (d) = {d + a_1 - 1 \choose a_1 - 1}
+
{d + a_2 - 2 \choose a_2 - 1}
+ \cdots +
{d + a_t - t \choose a_t - 1},
\end{equation}
where $t \leq s$ is largest such that $a_t \neq 0$.

There exists a projective $k$-scheme with Gotzmann representations (\ref{eq2}) and (\ref{eq4}), so $t=r$ and $b_i = a_i - 1$ for $i=1, \ldots, r$.

From the induction hypothesis, $\tilde{N'}$ is $r$-regular, so for all $q \geq 1$, we have 
\[
H^q (\tilde{N'}(d-q)) = 0
\]
for all $d \geq r$, so 
\[
H^q(\tilde{N}(d-q)) = 0
\]
for all $q \geq 2$ and all $d \geq s$. To show that this also holds for $q = 1$ and all $d = s$, suppose to the contrary that $H^1(\tilde{N}(s-1)) \neq 0$. 

Let $\mathfrak{m}$ be the maximal homogeneous ideal of $S$. Recall that the Hilbert  function and Hilbert polynomial can be related to the local cohomology:
\begin{equation}\label{eq-local-co}
H(M,d) - P_M(d) = \sum_{i}(-1)^i H_{\mathfrak{m}}^iM_d .
\end{equation}

Since $N$ is saturated and $F$ is $0$-regular, $H_{\mathfrak{m}}^0(M) = 0$ and $H_{\mathfrak{m}}^i(M) = H_{\mathfrak{m}}^{i+1}(N)$ for all $i \geq 1$. From the relations $H_{\mathfrak{m}}^i(N) = \bigoplus_d{H^{i-1}(\tilde{N}(d))}$ for $i \geq 2$, and equation \ref{eq-local-co}, $H(M,s-1) < P_M(s-1)$. So we can write
\[
H(M,s-1) \leq {s-1 + a_1 \choose s-1} + \cdots + {1+a_{s-1} \choose 1},
\]
and so 
\begin{align*}
H(M,s) & \leq H(M,s-1)^{\langle s-1 \rangle} \\
   & \leq {s + a_1 \choose a_1} + \cdots + {2+a_{s-1} \choose a_{s-1}} \\
   & < P_M(s).
\end{align*}

By repeating the last step, we have $H(M,d) < P_M(d)$ for all $d \geq s$, a contradiction. Therefore, $\tilde{N}$ is $s$-regular. $\qed$

The following corollary on regularity for arbitrary coherent sheaves is immediate.

\begin{corollary} \label{cor1} Suppose $\mathcal{F}$ is a coherent sheaf on $\mathbb{P}^n$ and $a \in \mathbb{Z}^{\geq 0}$ such that $\mathcal{F}(a)$ is generated by global sections and $\mathcal{F}(a)$ has Gotzmann regularity $s$. Then $\mathcal{F}$ is $s+a$-regular. $\qed$
\end{corollary}

\section{An explicit construction of the Quot scheme}

Let us use the above results and those of Gasharov to extend Gotzmann's construction of the Hilbert scheme to a construction of the Quot scheme. Write $\mathcal{O} = \mathcal{O}_{\mathbb{P}^n}$. For simplicity, consider $\mathrm{Quot}_P(\mathcal{O}^r)$, the quotients of $\mathcal{O}^r$ with Hilbert polynomial $P$. By Proposition \ref{prop1}, $P$ has a Gotzmann representation, say with Gotzmann number $s$. By Proposition \ref{prop2}, any coherent sheaf which is a quotient of $\mathcal{O}^r$ with Hilbert polynomial $P$ is generated in degrees at most $s$, and all cohomology vanishes in degrees at least $s$. 

The construction of the Quot scheme follows as in the construction of the Hilbert scheme. Let 
\begin{equation} \label{eqGrass}
\mathcal{G}_s = \mathrm{Gr}(S^r_s, P(s))
\end{equation}
 be the Grassmannian of $P(s)$-codimensional vector subspaces  of $S^r_s$, and similarly for $\mathrm{Gr}(S^r_{s+1}, P(s+1))$; and 
\begin{equation} \label{eqQuot}
W = \{ (F,G) \in\mathcal{G}_s \times \mathcal{G}_{s+1} | F \cdot S_1 = G\}.
\end{equation} 

\begin{theorem}\label{thmQuot}
The Quot scheme $\mathrm{Quot}_P(\mathcal{O}^r)$ is isomorphic to $W$.
\end{theorem}

The proof will follow \cite[Satz 3.4]{Got78}. Before beginning the proof, let us fix some notation. Write $\mathrm{Quot}_P^r$ for $\mathrm{Quot}_P(\mathcal{O}^r)$. Write $\mathbb{W}$ for the functor that will be represented by $W$, and  $\mathbb{G}_s$ for the functor of $\mathcal{G}_s$. For an affine scheme $T = \mathrm{Spec} A$ and the projection $\pi: T \times \mathbb{P}^n \to T$, recall that an element of $\mathbb{W}(A)$ is a pair $(F,G)$ of locally free direct summands of $S_s^r \otimes A$ and $S_{s+1}^r \otimes A$ of co-ranks $P(s)$ and $P(s+1)$ respectively, such that $F \cdot S_1 = G$.

Two lemmas will be given first. Their proofs are the same as the case of Hilbert schemes and ideals as given in their sources.

\begin{lemma}\label{lemmaQuot}
For any $d \geq s$ with $s$ being the Gotzmann number of $P$, the subfunctor $\mathbb{W}(A)$ of $\mathbb{G}_s(A)$ is representable by the scheme $W$, and the first projection embeds $W$ in $\mathcal{G}_s$; and $W$ contains the image of $\mathrm{Quot}_P^r$ under the closed diagonal embedding.
\end{lemma}

\proof The same argument from \cite[Proposition C.28]{Iar99} for the Hilbert scheme carries over to the Quot scheme. $\qed$

Thus we can use 
\[
W' = \{ F \in\mathcal{G}_s | \mathrm{corank}(F \cdot S_1) = P(s+1)\}
\]
interchangeably with $W$.

\begin{lemma}\label{lemmaQuot2}
If $A$ is a local Noetherian ring with maximal ideal $\mathfrak{m}$ and $A/\mathfrak{m} = k$, $S=A[x_0, \ldots, x_n]$, $N$ is a graded submodule of $S^r$ generated by $N_d$ such that $\mathrm{sat}(N \otimes k)$ is $d$-regular, and $M=F/N$ with $\dim M_d \otimes k = P(d)$, then the following are equivalent:
\begin{enumerate}
\item $M_n$ is free of rank $P(n)$ for all $n \geq d$
\item $M_n$ is free of rank $P(n)$ for all $n \gg 0$
\item $M \otimes k$ has Hilbert polynomial $P$, and $M_{d+1}$ is flat with rank $P(d+1)$.
\end{enumerate}
\end{lemma}

\proof The argument is the same as \cite[Satz 1.5]{Got78}. The equivalence of 1 and 2 is by \cite[Corollaire 7.9.9]{EGA3}, and clearly $1 \Rightarrow 3$.

To show $3 \Rightarrow 2$, assume that $M_{d+1}$ is flat. Writing $R$ for the syzygy module of $N$ (with an exact sequence $0 \to R \to S^m \to N \to 0$), one may argue as in Section 1.4 of \cite{Got78} to show that $R \otimes k$ is generated by $R_1 \otimes k$, or we may note that $\mathrm{sat}(N \otimes k)_{\geq d}$ has a linear resolution by  \cite[Proposition 1.3.1]{Cha07}, and that $\mathrm{sat}(N \otimes k)$ and $N \otimes k$ agree in degree $d$ and larger, since they have the same Hilbert function for $n \geq d$. 

Now if $v \gg 0$ and $x \otimes 1 \in \ker(N_{d+v} \otimes k \to S^r_{d+v} \otimes k)$, then writing $x=fy$ with $f \in N_d$ and $y \in S^m_v$, we have $y \otimes 1 \in R_v \otimes k = (S_{v-1} \otimes k)R_1$, hence $y \in S_{v-1}R_1 + \mathfrak{m}S^m_v$. It follows that $x \in \mathfrak{m} N_{d+v}$, thus $N_{d+v} \otimes k \to S^r_{d+v} \otimes k$ is an injection; therefore a basis for $M_{d+v} \otimes k$ can be lifted to a basis for $M_{d+v}$ and $M_{d+v}$ is free. $\qed$

\proof (of Theorem \ref{thmQuot})

Let $(F,G)$ be the universal element of $\mathbb{W}(W)$ and set $M = \oplus_d (S_{s+d}^r/S_d \cdot F)$. Note that $\mathrm{Quot}_P^r$ is the maximal subscheme of $W$ such that $M \otimes_W \mathcal{O}_{\mathrm{Quot}_P^r}$ is flat over $\mathrm{Quot}_P^r$ with Hilbert polynomial $P$. For each $d \geq s$, take $Z_d$ to be the maximal subscheme of $W$ such that $M_d \otimes \mathcal{O}_{Z_d}$ is flat over $\mathcal{O}_{Z_d}$ with rank $P(d)$. By Theorem \ref{Gas97thm} \#1, $\dim M_d \otimes k(p) \leq P(d)$ for all $p \in W$ and $d \geq s$, hence $Z_d$ is closed. There exists $d_0$ such that $Z_{d_0} = Z_d$ for all $d \geq d_0$, since the $Z_d$ form an ascending chain, so $\mathrm{Quot}_P^r = Z_{d_0}$ is closed.  

Now consider the set $Z$ of points $p \in W$ such that $M \otimes_W k(p)$ has Hilbert polynomial $P$. Note that $Z_{d_0} \subseteq Z$. By Lemma \ref{lemmaQuot2}, for any $p \in Z$, we have $M_d \otimes_W \mathcal{O}_p$ is free of rank $P(d)$ over $\mathcal{O}_p$ for all $d \geq s$, so $M \otimes_W \mathcal{O}_p$ is free over $\mathcal{O}_p$ with Hilbert polynomial $P$, thus $p \in Z_{d_0}$. It follows that $Z$ is the set of points for which $M_{d_0} \otimes_W \mathcal{O}_p$ is free of rank $P(d_0)$ over $\mathcal{O}_p$, so $Z$ is an open subset of $W$, and thus $\mathrm{Quot}_P^r$ is an open subscheme of $W$. 

Since $\mathrm{Quot}_P^r$ is an open and closed subscheme of $W$, it remains to show that all connected components of $W$ are contained in $\mathrm{Quot}_P^r$, for which it is sufficient to check geometric points. Let $(F,G) \in \mathbb{W}(K)$ be a geometric point. Then $M = \oplus_d (S_{s+d}^r/S_d \cdot F)$ has Hilbert function $P(s)$ and $P(s+1)$ in degrees $s$ and $s+1$, so by Theorem \ref{Gas97thm} \#2, $M$ has Hilbert polynomial $P$, so $(F,G)$ is in the image of $\mathrm{Quot}_P^r$. Therefore, $\mathrm{Quot}_P^r \cong W$.
$\qed$

\section{Quot schemes on $\mathbb{P}^1$}

Consider this construction for the Quot scheme
\[
\mathcal{Q} = \mathrm{Quot}_P(\mathcal{O}_{\mathbb{P}^1}^r),
\]
with $P(d) = k(d+1) + m$. The Gotzmann representation of $P$ is
\[
{d+1 \choose 1} + {d \choose 1} + \cdots + {d-k+2 \choose 1} + \frac{k(k-1)}{2} + m,
\]
so the Gotzmann number of $P$ is 
\[
s = \frac{k(k+1)}{2} + m.
\]
Thus, 
\[
\mathcal{Q} \cong \{F \in \mathrm{Gr}(r(s+1), k(s+1) + m) | \mathrm{codim}F \cdot S_1 = k(s+2)+m\}.
\]

Let $\mathcal{F} \in \mathcal{Q}$, and let $l$ be an integer at least the Castelnuovo-Mumford regularity of $\mathcal{F}$. We have an exact sequence
\[
0 \to \mathcal{K}(l) \to \mathcal{O}^r(l) \to \mathcal{F}(l) \to 0.
\]
Set $V = H^0(\mathcal{K}(l)) \subset S_l^r$. Then $\dim V = (r-k)(l+1) - m$. If $\dim V \cdot S_1 = (r-k)(l+2) - m$, one expects an element $F$ in $\mathrm{Gr}((r-k)(l+1)-m, r(l+1))$ and an element $G = F \cdot S_1$ in $\mathrm{Gr}((r-k)(l+2)-m, r(l+2))$. The Porteus formula gives a degeneracy locus with codimension  
\begin{align*}
\ & [2(r-k)(l+1)-2m - (r-k)(l+2) + m][r(l+2) - (r-k)(l+2)+m] \\
  = & [(r-k)l - m][k(l+2) + m].
\end{align*}
So the expected dimension of $\mathcal{Q}$ is
\[
[(r-k)(l+1)-m][k(l+1)+m] - [(r-k)l-m][k(l+2) + m] = (r-k)k + rm.
\]
On the other hand, given
\[
0 \to \mathcal{K} \to \mathcal{O}^r \to \mathcal{F} \to 0
\]
we can write $\mathcal{K} = \oplus \mathcal{O}(-t_i)$, and since the Hilbert polynomial of $\mathcal{F}$ is $k(d+1) +m$, we have $r-k$ summands of $\mathcal{K}$, each $t_i$ is non-negative, and $\sum{t_i} = m$. In order to compute $\dim Hom(\oplus \mathcal{O}(-t_i),\mathcal{O}^r)/Aut(\oplus \mathcal{O}(-t_i))$, we need a bound on $\dim Aut(\oplus \mathcal{O}(-t_i))$.

Rewrite $\mathcal{K} = \oplus \mathcal{O}(-t_i)$ by gathering common twists, as
\[
\mathcal{K} = \oplus_{i=0}^s{\mathcal{O}(i)^{e_i}}.
\]
Then $\dim Aut{K} = \sum_{0 \leq i \leq j \leq s}{(j-i+1)e_ie_j}$. This dimension is bounded as follows.

\begin{lemma}
Assume $m$ and $n$ are integers such that $0 \leq m \leq n$ and $e_0, e_1, \ldots, e_n$ are such that $\sum{e_i} = m$ and $\sum{i e_i} = n$. Then 
\[
\sum_{0 \leq i \leq j \leq n}{(j-i+1)e_ie_j} \geq m^2,
\]
with a unique solution giving equality.
\end{lemma}
\proof Observe that 
\begin{align*}
\sum_{0 \leq i \leq j \leq n}{(j-i+1)e_ie_j} & = \sum_{i=0}^n{e_i^2} + \sum_{0 \leq i < j \leq n}{(j-i+1)e_ie_j} \\
  & \geq \sum_{i=0}^n{e_i^2} + \sum_{0 \leq i < j \leq n}{2e_ie_j} \\
	& = m^2,
\end{align*}
and equality occurs only if there are at most two non-zero exponents. Assume $0 < e_i \leq m$, and $e_i + e_{i+1} = m$ and $i e_i + (i+1) e_{i+1} = n$. Then $n = i(e_i + e_{i+1}) + e_{i+1} = im + e_{i+1}$, with $0 \leq e_{i+1} < m$. By the Division Algorithm, such a solution exists and is unique. $\qed$

So 
\begin{align*}
\dim & Hom(\oplus \mathcal{O}(-t_i),\mathcal{O}^r)/Aut(\oplus \mathcal{O}(-t_i)) \\
  & \  = \dim H^0(\oplus \mathcal{O}(e_i)^r) - \dim Aut(\oplus \mathcal{O}(-t_i)) \\
  & \ \leq r(r-k) - rm - (r-k)^2 \\
	& \ = (r-k)k + rm,
\end{align*}
agreeing with the previous computation; and the uniqueness makes $\mathcal{Q}$ irreducible of the expected dimension. The result of this computation gives the following:

\begin{prop}The Quot Scheme $\mathrm{Quot}_P(\mathcal{O}_{\mathbb{P}^1}^r)$ with $P(d) = k(d+1) + m$ having Gotzmann number
\[
\frac{r(r+1)}{2} + m,
\]
is the irreducible scheme given by the degeneracy locus of 
\[
\mathrm{Gr}((r-k)(\frac{r(r+1)}{2}+m+1)-m,r(\frac{r(r+1)}{2}+m+1))
\]
with rank $(r-k)k + rm$. $\qed$
\end{prop}

Iarrobino and Kleiman \cite{Iar99} showed that the maximal Castelnuovo-Mumford regularity of a scheme with Hilbert polynomial $P$ is equal to the Gotzmann number of $P$. However, this is not the case for globally generated sheaves. In fact, the following example shows that we can use the Gotzmann construction to embed a Quot scheme into a smaller Grassmannian than the Gotzmann number gives. 

\begin{ex}
Consider $\mathcal{Q} = \mathrm{Quot}_P(\mathcal{O}_{\mathbb{P}^1}^3)$, where $P(d) = 2(d+1)$. The Gotzmann number of $P$ is 3, and the Castelnuovo-Mumford regularity of a sheaf $\mathcal{F} \in \mathcal{Q}$ is 0. Set
\begin{align*}
W_0 & = \{ F \in \mathrm{Gr}(S_0^3,P(0)) | \codim F \cdot S_1 = P(1)\} \\
   & = \{ F \in \mathrm{Gr}(3,2) | \codim F \cdot S_1 = 4\} \\
	 & \cong \mathbb{P}^2
\end{align*}
and
\begin{align*}
W_3 & = \{ F \in \mathrm{Gr}(S_3^3,P(3)) | \codim F \cdot S_1 = P(4)\} \\
   & = \{ F \in \mathrm{Gr}(12,8) | \codim F \cdot S_1 = 10\}.
\end{align*}
Note that $\mathcal{Q} \cong W_3$, but that by the last section, $\mathcal{Q}$ embeds into $W_0$. Given a point $F \in W_0$, one shows that $F \cdot S_3 \in W_3$, so $Q \cong W_0$. 
\end{ex}

\section{Bounding Chern classes of globally generated coherent sheaves}

Given a rank 2 globally generated vector bundle $\mathcal{F}$ with Chern classes $c_1$ and $c_2$, \cite{Chi12} gives a bound for $c_2$ in terms of $c_1$, using the vanishing of a section of $\mathcal{F}$ along a smooth curve. Recall that a vector bundle is \emph{split} if it can be written as a direct sum of line bundles.

\begin{prop}(\cite{Chi12}, Lemma 1.5)\label{propChi12}
If $\mathcal{E}$ is a non-split rank-2 globally generated vector bundle on $\mathbb{P}^3$ with $c_1 \geq 4$, then 
\[
c_2 \leq \frac{2c^3_1 - 4c^2_1 + 2}{3c_1 - 4}. \qed
\]
\end{prop}

We can use the results herein to give a larger bound, but a bound that applies to any globally generated coherent sheaf.

Let $\mathcal{E}$ be a rank 2 globally generated coherent sheaf on $\mathbb{P}^3$ with Chern classes $c_1$, $c_2$ and $c_3$. Note that $c_3$ may be non-zero, since we are not restricting attention to vector bundles. We compute the Hilbert polynomial using Hirzebruch-Riemann-Roch:

\begin{equation}\label{eqHilb-Chern}
P_{\mathcal{E}}(d) = \frac{1}{3}d^3 + (2 + \frac{1}{2}c_1)d^2 + (\frac{1}{2}c_1^2 + 2c_1 + \frac{11}{3} - c_2)d + \overline{c},
\end{equation}
where 
\[
\overline{c} = c_1^3/6 + c_1^2 + 11c_1/6 - c_1c_2/2 - 2c_2 + c_3/2 + 2
\]
is given for completeness, although it has no bearing on the following computation.

Since $\mathcal{E}$ is globally generated, $P_{\mathcal{E}}(d)$ has a Gotzmann representation; write

\[
P_{\mathcal{E}}(d) = P_3 + P_2 + P_1 + P_0,
\]
where $P_i$ is the sum of binomial coefficients in the Gotzmann representation of degree $i$. Since $\mathcal{E}$ is rank $2$, it follows that

\[
P_3 = {d + 3 \choose 3} + {d + 2 \choose 3} = \frac{1}{3}d^3 + \frac{3}{2}d^2 + \frac{13}{6}d + 1,
\]
so 

\[
P_{\mathcal{E}} - P_3 = \frac{1}{2}(1 + c_1)d^2 + (\frac{1}{2}c_1^2 + 2c_1 + \frac{3}{2} - c_2)d + \overline{c} - 1.
\]

This tells us that $c_1 \geq -1$, but in fact $c_1 \geq 0$ for any globally generated coherent sheaf on $\mathbb{P}^n$. To see this, assume $\mathcal{F}$ is a rank $r$ globally generated sheaf on $\mathbb{P}^n$ with Hilbert polynomial $P(d)$; the top two terms of $P(d)$ are $rd^n/n! + (r(n+1)/2 + c_1(\mathcal{F}))d^{n-1}/(n-1)!$. By reason of rank, the torsion subsheaf $\mathrm{tors}(\mathcal{F})$ has a Hilbert polynomial of degree at most $n-1$; so the Hilbert polynomial of the torsion-free quotient $\overline{\mathcal{F}}$ shows that $c_1(\mathcal{F}) \geq c_1(\overline{\mathcal{F}})$. 

So now assume in addition that $\mathcal{F}$ is torsion-free. For a hyperplane $i:H \hookrightarrow \mathbb{P}^n$, there is an exact sequence
\[
0 \to \mathcal{F}(-1) \to \mathcal{F} \to \mathcal{F}_H \to 0,
\]
and $c_1(\mathcal{F}) = c_1(\mathcal{F}_H)$. Since there exists a line $L$ such that $\mathcal{F}_L$ is a vector bundle, the result $c_1 \geq 0$ follows from induction on $n$.

Next, it follows that $P_2$ has $c_1 + 1$ terms:

\[
P_2 = {d \choose 2} + \cdots + {d - c_1 \choose 2} = \frac{1}{2}(1 + c_1)d^2 - \frac{1}{2}(c_1 + 1)^2 d + b,
\]
where $b = \frac{1}{6}(c_1^3 + 3c_1^2 + 2c_1)$. So

\[
P_{\mathcal{E}} - P_3 - P_2 = (c_1^2 + 3c_1 + 2 - c_2)d + \overline{c} - 1 - b.
\]

This leading coefficient must be non-negative for there to be a Gotzmann representation. Therefore we have shown:

\begin{prop}\label{propChern} Let $\mathcal{E}$ be a rank 2 globally generated coherent sheaf on $\mathbb{P}^3$. Then its first and second Chern classes satisfy the inequality 
\[
c_2 \leq c_1^2 + 3c_1 + 2. \qed
\]
\end{prop}

\begin{ex}
Let $S=\mathbb{C}[x,y,z,w]$, and take $N$ to be the image of the matrix
\[\left(
\begin{array}{lll}
x^4 & y^4 & z^4 \\
x^4 & y^4 & 0
\end{array}
\right) 
\]
and let $\mathcal{E}$ be the sheaf associated to $N$. Then $\mathcal{E}(4)$ is globally generated, and from the resolution
\[
0 \to S(-4) \to S^3 \to N(4) \to 0
\]
it follows that $\mathcal{E}(4)$ has Hilbert polynomial
\[
P_{\mathcal{E}(4)}(d) = \frac{1}{3}d^3 + 4d^2 + \frac{11}{3}d + 4.
\]
From equation \ref{eqHilb-Chern}, one computes the Chern classes $c_1(\mathcal{E}(4)) = 4$, $c_2(\mathcal{E}(4)) = 16$ and $c_3(\mathcal{E}(4)) = 64$. Note that $\mathcal{E}(4)$ is not a vector bundle, and in fact, it does not satisfy the inequality for vector bundles in Proposition \ref{propChi12}; but it does satisfy the inequality in Proposition \ref{propChern}.
\end{ex}

\section*{Acknowledgements}

The author wishes to thank Vladimir Baranovsky for helpful discussions on this work, and the anonymous referee for important corrections and suggestions.

\bibliographystyle{elsarticle-num}
\bibliography{Gotzmann_rep}

\end{document}